\documentclass[11pt]{amsart}
\usepackage{epsfig}
\usepackage{xy}
\usepackage{graphicx,color}
\usepackage{amsmath,amssymb}
\usepackage{float}
\usepackage{graphics}

\xyoption{all}

\newtheorem{lemma}{Lemma}[section]
\newtheorem{theorem}[lemma]{Theorem}

\newcommand{\Hom}{\operatorname{Hom}}

\newcommand{\Ext}{\operatorname{Ext}}
\newcommand{\Tor}{\operatorname{Tor}}

\newcommand{\End}{\operatorname{End}}

\newcommand{\add}{\operatorname{add}}
\newcommand{\module}{\operatorname{mod}}

\renewcommand{\P}{\mathcal P}

\newcommand{\T}{\mathcal T}

\newcommand{\C}{\mathcal C}
\newcommand{\D}{\mathcal D}

\newcommand{\isom}{\cong}

\title{An introduction to Higher cluster categories}

\author[Buan]{Aslak Bakke Buan}
\address{Institutt for matematiske fag\\
Norges teknisk-naturvitenskapelige universitet\\
N-7491 Trondheim\\
Norway}
\email{aslakb@math.ntnu.no}
\thanks{The author is supported by an NFR-frinat grant}
\begin{document}

\maketitle

\section*{Introduction}

Cluster categories were defined in \cite{bmrrt} in order to use
categorical methods to give a conceptual model for the combinatorics
of cluster algebras, as defined by Fomin and Zelevinsky
\cite{fz}. 
With contributions from many mathematicians, this theory and its
generalisations have given new links between
categorical representation theory and several branches of mathematics
and mathematical physics. In addition, various problems 
concerning cluster algebras and related combinatorial problems have
been solved. 
There are several recent survey papers on this
topic, e.g. \cite{k2, k3, rei}, discussing both categorical and combinatorial aspects of
the theory. 
 
In this survey we discuss some combinatorial aspects of 
a generalisation of cluster categories, called $m$-cluster
categories, or higher cluster categories. Such categories are not
explicitly linked to cluster algebras.
A survey on categorical aspects of  higher cluster categories, and
generalisations, is given in \cite{k4}.

A cluster category is defined as an orbit category
of the derived category of an hereditary finite dimensional
algebra. Loosely speaking, it is obtained by identifying the AR-translation
$\tau$ with the shift $[1]$.
Keller \cite{k} proved that a cluster
category is triangulated, and that the canonical functor from the
derived category to the cluster category is a triangle functor.
The orbit category is a Calabi-Yau category of CY-dimension 2.

Keller's proof also showed that the orbit category obtained by identifying
$\tau$ with the $m$-fold shift $[m]$ is triangulated. 
These categories has later been called $m$-cluster categories, and
they are Calabi-Yau of dimension $m+1$.

The main interest in 
$1$-cluster categories, and some other triangulated categories of CY-dimension 2, is due to the
combinatorial properties of the set of tilting objects (also called cluster
tilting objects).  
The definition of tilting objects canonically extends to $m$-cluster categories.

In this survey, we give an overview over some combinatorial
aspects of the set of tilting objects in an $m$-cluster category, with
focus on those properties which are valid for all $m\geq 1$.
In the two first sections
we give some more details and background and a
precise definition. We also recall definitions and results on tilting
theory in higher cluster categories. The results in these sections are
mainly due to Wraalsen, Zhou and Zhu \cite{w,z,zz}.
Then, in the next three sections,
we consider three different, but related, combinatorial aspects of
the set of tilting objects $\T$ in $m$-cluster categories. 
First, we discuss work of Baur and Marsh, who model the combinatorics
of $\T$ in
the Dynkin case $A$ or $D$ using arcs in certain (unpunctured or
punctured) polygons \cite{bm1,bm2}.
Next, we discuss links to the Fomin-Reading generalised
associahedra \cite{fr}, due to Thomas \cite{t} and Zhu \cite{z}.
Then, in section 5, we explain coloured quivers and mutation of such, as defined in
joint work with Thomas \cite{bt}, and show how these can be used to
describe combinatorial aspects of $\T$ for arbitrary finite quivers.
We end, in section 6, with some comments on other aspects of higher
cluster categories and generalisations.

The author would like to thank the organisers of
the CIMPA-UNESCO-IPM School in 
Representation Theory of Algebras in June 2008 in Tehran, and especially
Javad Asadollahi and his colleagues at IPM for their warm hospitality.
He would also like to thank Hermund Torkildsen for commenting on an 
earlier version of this
survey, and to an anonymous referee for several good suggestions concerning 
the presentation.

\section{Background and definition}

We give some background on derived categories, before we discuss the
construction of the $m$-cluster categories. For more information on
derived categories, see \cite{h, tria}. For basic information on
finite dimensional algebras and their representation theory, see the textbooks \cite{ars,ass}. 

\subsection{The derived category}

Let $H$ be a hereditary finite dimensional algebra over an
algebraically closed field $k$. We assume $H$ is basic, hence
$H$ is isomorphic to a path algebra $kQ$ of some finite quiver $Q$.
Let  $\module H$ be the category of finite dimensional left
$H$-modules, and
let $D^b(H)$ be the (bounded) derived category.
Let $[1]$ denote the shift functor on $D^b(H)$, let $[-1]$ denote its
inverse.
The derived category is a 
Krull-Schmidt category, and its indecomposable objects are isomorphic
to stalk complexes $M[i]$, where $M$ is an indecomposable $H$-module,
and $i$ is some integer.
For indecomposables $M[i]$ and $N[j]$, we have that the morphism spaces
are given by
$$\Hom_{D^b(H)}(M[i], N[j]) = 
\begin{cases}
\Hom_H(M,N) & \text{if } $i=j$ \\
\Ext^1_H(M,N) & \text{if } $j= i+1$  \\
0 & \text{else.}
 \end{cases}
$$

By results of Happel \cite{h}, the derived category  $D^b(H)$  has
Auslander-Reiten triangles. This implies that there is an
autoequivalence $\tau$ on the derived category, with the property that
for each  indecomposable object $M$, there is a uniquely determined
triangle
$$ \tau M \to E \to M \to.$$
Furthermore, we have the Auslander-Reiten
formula $$\Hom_{D^b(H)}(M,N[1]) \simeq D\Hom_{D^b(H)}(N, \tau M),$$
where $D=\Hom(\ ,k)$ is the ordinary duality.

We view objects in $\module H$ as stalk complexes in degree 0.
If $M$ is a non-projective indecomposable module, then $\tau M$
coincides with $\tau_H M$, where $\tau_H$ denotes the AR-translation in
the module category.
If $P$ is an indecomposable projective, then $\tau P = I [-1]$, where
$I = D\Hom_H(P,H)$ is indecomposable and injective. 

\subsection{An example of type $A$}\label{ex1}
Let $Q$ be the quiver
$$
\xymatrix{
1 \ar[r] & 2 & 3 \ar[l] \ar[r]& 4
}
$$

Consider the path algebra $H =kQ$, and let $e_i$ be the idempotent in $H$
corresponding to the vertex $i$. There are 10 indecomposable modules
in $\module H$. These are the 4 projectives $P_i = He_i$ and the
4 injectives $I_i = D(e_iH)$, in addition to the two modules $X \isom
(P_1 \amalg P_3)/P_2$  and $Y= P_3/P_4$. The AR-quiver of the module category is
given as follows, where the action of $\tau$ is indicated by the dotted arrows.

\def \objectstyle{\scriptstyle}
\def \labelstyle{\scriptstyle}
$$
 \xymatrix@!0@C1.20cm@R0.7cm{ 
&P_1 \ar[dr] & & I_4 \ar[dr] \ar@{.>}[ll] & \\
P_2 \ar[dr] \ar[ur]  & & X \ar[ur] \ar[dr] \ar@{.>}[ll] & & I_3 \ar@{.>}[ll]   \\
& P_3 \ar[ur] \ar[dr] & & I_2 \ar@{.>}[ll] \ar[dr] \ar[ur]  &  \\
P_4  \ar[ur]  & & Y \ar[ur] \ar@{.>}[ll] & & I_1 \ar@{.>}[ll]
}
$$

In the derived category, the AR-translation $\tau$ is defined on all objects, and actually
becomes an autoequivalence.
A segment of the AR-quiver of the derived category looks as follows,
where for an indecomposable $M$,
we have that $\tau M$ is the neighbour directly to left.

\def \objectstyle{\scriptstyle}
\def \labelstyle{\scriptstyle}
$$
\xymatrix@!0@C0.83cm@R0.7cm{ 
& & I_1[-1] \ar[dr] & & P_1 \ar[dr] & & I_4 \ar[dr] & & P_4[1] \ar[dr]
& & Y[1] \ar[dr]  & &  I_1[1] \ar[dr]
& & P_1[2] \ar[dr] & & P_4[2] \\
\cdots & I_2[-1]  \ar[dr] \ar[ur] & & P_2 \ar[dr] \ar[ur]  & & X
\ar[dr] \ar[ur]  & & I_3 \ar[dr] \ar[ur]   & & P_3[1] \ar[dr]
\ar[ur]   & & I_2[1] \ar[dr] \ar[ur]  & &
P_2[2]  \ar[dr] \ar[ur]  & & X[2] \ar[dr] \ar[ur]   & \cdots \\
& & I_3[-1] \ar[ur]  \ar[dr]  & & P_3 \ar[ur]  \ar[dr]  & & I_2 \ar[ur]  \ar[dr]  & & P_2[1]
\ar[ur]  \ar[dr]   &  & X[1] \ar[ur]  \ar[dr]  & &  I_3[1] \ar[ur]  \ar[dr]  
& & P_3[2] \ar[ur]  \ar[dr]  & & I_2[2] \\
\cdots & I_4[-1]  \ar[ur] & & P_4  \ar[ur]  & & Y
\ar[ur]  & & I_1  \ar[ur]   & & P_1[1] 
\ar[ur]   & & I_4[1]  \ar[ur]  & &
P_4[2]  \ar[ur]  & & Y[2]  \ar[ur]   & \cdots 
}
$$

\subsection{The $m$-cluster category}

Consider now the autoequivalence $G = \tau^{-1}[m]$ on $D^b(H)$, and define the
$m$-cluster category to be the orbit category $\C = D^b(H)/G$. 

The objects of $\C$ are the $G$-orbits of
objects in $D^b(H)$; we use the same notation for an object in
$D^b(H)$ and its orbit in $\C$. The morphism spaces in $\C$,
are given by 
$$\Hom_{\C} (X,Y)= \amalg_i \Hom_{\D^b(H)}(X,  G^i Y) $$

Keller \cite{k} proved that $\C$ is triangulated, and that the
canonical functor $D^b(H) \to \C$ is a triangle functor. 
It follows from \cite{bmrrt} that $\C$ is a Krull-Schmidt category
with almost split triangles and translation functor induced from
$D^b(H)$, and it can be shown that the AR-formula
$$\Hom_{\C}(M,N[1]) \simeq D\Hom_{\C}(N, \tau M),$$
still holds in $\C$.

There is a canonical embedding of $\module H$ into $D^b(H)$.
Let $\module H[0]$ denote the image under this embedding, and let 
$\module H[i]$ be defined in the obvious way. 
We say that $\module H[0] \vee \dots \vee \module H[m-1] \vee H[m]$ is a
standard domain in $D^b(H)$. It is clear from the definition of $\C$,
that any indecomposable object in $\C$ is up to isomorphism induced by an object in the
standard domain. 

\subsection{Example}\label{arex}
We consider the path algebra of example \ref{ex1}. 
Now the $2$-cluster category is of finite type, consisting of
24 indecomposable objects: 2 copies of the 10 indecomposable objects in the module category and one additional copy of
the 4 indecomposable projectives. The AR-quiver
looks as follows, where
one should note that objects on the left border and the right border are identified.
\def \objectstyle{\scriptstyle}
\def \labelstyle{\scriptstyle}
$$
\xymatrix@!0@C0.83cm@R0.7cm{ 
 & P_1 \ar[dr] & & I_4 \ar[dr] & & P_4[1] \ar[dr]
& & Y[1] \ar[dr]  & &  I_1[1] \ar[dr]
& & P_1[2] \ar[dr] & & P_1 \\
 P_2 \ar[dr] \ar[ur]  & & X
\ar[dr] \ar[ur]  & & I_3\ar[dr] \ar[ur]   & & P_3[1] \ar[dr]
\ar[ur]   & & I_2[1] \ar[dr] \ar[ur]  & &
P_2[2]  \ar[dr] \ar[ur]  & & P_2 \ar[dr] \ar[ur]   &  \\
& P_3 \ar[ur]  \ar[dr]  & & I_2 \ar[ur]  \ar[dr]  & & P_2[1]
\ar[ur]  \ar[dr]   &  & X[1] \ar[ur]  \ar[dr]  & &  I_3[1] \ar[ur]  \ar[dr]  
& & P_3[2] \ar[ur]  \ar[dr]  & & P_3 \\
 P_4  \ar[ur]  & & Y
\ar[ur]  & & I_1  \ar[ur]   & & P_1[1] 
\ar[ur]   & & I_4[1]  \ar[ur]  & &
P_4[2]  \ar[ur]  & & P_4  \ar[ur]   &  
}
$$

\section{Tilting objects and exchange triangles}\label{tilt}

Tilting theory in module categories over
finite dimensional algebras was initiated more than 30 years ago, see
\cite{handbook}.
The original motivation was to compare module categories. Happel \cite{h}
introduced the use of derived categories in the theory, and showed
that algebras related by tilting are derived equivalent.

In the setting of hereditary algebras, a tilting module in $\module H$ is a module $T$ with $\Ext_{H}^1(T,T) =
0$ and with $n$ indecomposable non-isomorphic direct summands,
where $H$ has $n$ isomorphism-classes of simples. 

In work of Riedtmann and Schofield \cite{rs}, Unger \cite{u}, and others,
combinatorial properties on the set of tilting
modules were studied, in particular the simplicial complex defined by
the set of direct
summands in tilting modules was introduced. 
See \cite{handbook-unger} for more background on combinatorial aspects 
of tilting modules for finite dimensional algebras.

In this section we will define tilting objects in (higher) cluster category.
Using the natural embedding of a module category into a cluster
category, it is easy to see that tilting modules will be mapped to
tilting objects. In fact, for $1$-cluster categories, all tilting
objects are of this form (up to
derived equivalence). In the case of higher cluster categories
there are more tilting objects, as we will observe in later examples.

\subsection{Tilting theory in $m$-cluster categories}

An object $M$ in an $m$-cluster category is called rigid if
$\Ext^i_{\C}(M,M) = 0$ for $i=1, \dots, m$. A finite collection of
rigid objects $\{X_i \}$
is said to be  $\Ext$-compatible if the direct sum $\amalg  X_i$ is
rigid.
$M$ is called maximal rigid if the indecomposable direct summands in $M$ form a
maximal $\Ext$-compatible collection.
A tilting object $M$ in $\C$ is a rigid object with the additional
property that if an object $X$ satisfies $\Ext^i(M,X)= 0$ for $i=1,
\dots, m$, then this implies that $X$ is in
$\add M$.

Zhu \cite{z}, see also \cite{w}, showed that tilting and maximal rigid
objects coincide. This was shown in \cite{bmrrt} in the case $m=1$.
Recall that an object $X$ is called {\em basic} if any indecomposable object occurs
at most once in a direct sum decomposition of $X$.

\begin{theorem}
\cite{z} The following are equivalent for a basic rigid object $T$ in an
$m$-cluster category.
\begin{itemize}
\item[(a)] $T$ is maximal rigid.
\item[(b)] $T$ is tilting.
\item[(c)] $T$ has $n$ isomorphism classes of indecomposable direct summands. 
\end{itemize}
\end{theorem}

Note that 
it follows from this that every (basic) rigid object is a direct summand in a
tilting object.

\subsection{Complements}

Let $T = \amalg_{i=1}^n  T_i $ be a tilting object in an $m$-cluster
category, and fix an indecomposable direct summand $T_k$. 

We call $B_k = T/T_k$ an {\em almost complete tilting object}, and
indecomposable objects $X$ such that $B_k \amalg X$ is tilting, are
called complements to $B_k$. Indeed,
$T_k$ is a complement. Let $T_k \overset{f}{\to} B_k'$ be a minimal left $\add
B_k$-approximation of $T_k$. This means:
\begin{itemize}
\item[-] $B_k'$ is in $\add B_k$.
\item[-] Any map from $T_k$ to an object in $\add B_k$, factors through
  the map $f$.
\item[-] If $gf = f$ for some endomorphism $g \colon B'_k \to B'_k$, then
  $g$ is an automorphism. 
\end{itemize}
Let 
\begin{equation}\label{exc}
T_k \to B'_k \to T_k^{\ast} \to
\end{equation}
be the induced triangle in $\C$. 
Then one can show that $T_k^{\ast}$ is also a complement to $B_k$ with
$T_k^{\ast} \not \simeq T_k$. 
The triangle (\ref{exc}) is called an
{\em exchange triangle}. One can of course iterate this
procedure to produce new complements  and exchange triangles. However,
one can show that after $m$ iterations, such that totally $m+1$ complements
are constructed, no new complements will occur. 
Also, one can show that $B_k$ has no further complements than those
constructed in this way.
More precisely, we have the following theorem.

\begin{theorem}
\cite{w,zz} The almost complete tilting object $B_k$ in the $m$-cluster category $\C$
has exactly $m+1$ complements $T_k^{(c)}$ for $c= 0,1, \dots , m$
occurring in $m$ exchange triangles
\begin{equation}\label{exch}  T_k^{(c)} \overset{f_k^{(c)}}{\rightarrow} B_k^{(c)}
\overset{g_k^{(c+1)}}{\rightarrow} T_k^{(c+1)} \overset{h_k^{(c+1)}}{\rightarrow}  \end{equation}
\end{theorem}

The fact that we get $m+1$ complements in this way was proved in
\cite{iy}, while the fact that there are no further complements was proved
independently in \cite{zz} and in \cite{w}. 

It is pointed out in \cite{zz}, that exchange is transitive on the set of
tilting objects; i.e any tilting object can be reached from any other
tilting object by a finite sequence of exchanges. This was proved in
\cite{bmrrt} for $m=1$, using ideas of \cite{happelunger}.

\subsection{Example}\label{exA}
We revisit our example \ref{ex1}.
The boxed objects are  the direct summands of an almost complete tilting
object $B = I_4 \amalg I_1 \amalg Y[1]$, and 
the encircled object are the three complements of $B$.

\def \objectstyle{\scriptstyle}
\def \labelstyle{\scriptstyle}
$$
\xymatrix@!0@C0.83cm@R0.7cm{ 
 & P_1 \ar[dr] & & *+[F]{I_4} \ar[dr] & & P_4[1] \ar[dr]
& & *+[F]{Y[1]} \ar[dr]  & &  I_1[1] \ar[dr]
& & P_1[2] \ar[dr] & & P_1 \\
 P_2 \ar[dr] \ar[ur]  & & *+[o][F]{X}
\ar[dr] \ar[ur]  & & I_3 \ar[dr] \ar[ur]   & & P_3[1] \ar[dr]
\ar[ur]   & & I_2[1] \ar[dr] \ar[ur]  & &
P_2[2]  \ar[dr] \ar[ur]  & & P_2 \ar[dr] \ar[ur]   &  \\
& P_3 \ar[ur]  \ar[dr]  & & I_2 \ar[ur]  \ar[dr]  & &  *+[o][F]{P_2[1]}
\ar[ur]  \ar[dr]   &  & X[1] \ar[ur]  \ar[dr]  & &   *+[o][F]{I_3[1]} \ar[ur]  \ar[dr]  
& & P_3[2] \ar[ur]  \ar[dr]  & & P_3 \\
 P_4  \ar[ur]  & & Y
\ar[ur]  & & *+[F]{I_1}  \ar[ur]   & & P_1[1] 
\ar[ur]   & & I_4[1]  \ar[ur]  & &
P_4[2]  \ar[ur]  & & P_4  \ar[ur]   &  
}
$$

\bigskip
The three exchange triangles are:
$$X \to I_1 \amalg I_4 \to P_2[1] \to $$
\ 
$$P_2[1] \to Y[1] \to I_3[1] \to $$
\
$$I_3[1] \to 0 \to I_3[2] (=X) \to$$

\section{A graphical description}

Independent of the ideas in \cite{bmrrt},
Caldero, Chapoton and Schiffler  \cite{ccs} defined a family of categories, using
diagonals in regular $n$-gons as objects. They also showed that their categories are
equivalent to the cluster categories of Dynkin type $A$. Later
Schiffler \cite{s} used a similar approach to describe the cluster
categories of type $D$. He considered punctured $n$-gons instead.

Generalising this, Baur and Marsh gave a graphical interpretation
of $m$-cluster categories in type $A$ \cite{bm2} and in type $D$
\cite{bm1}.  See \cite{baur} for a survey. Here we will give a brief discussion of their ideas in type $A$,
including an example. 

\subsection{A category from polygons}

We discuss here the results of Baur and Marsh \cite{bm2} for Dynkin type $A$.
We want to construct a certain category of diagonals of an 
$(nm+2)$-gon $\P = \P_{nm+2}$,   
where $m$ and $n$ are positive integers, and $n>1$. This category will be
equivalent to the $m$-cluster category of a Dynkin quiver of type $A_{n-1}$.
The indecomposable objects in the $m$-cluster
category $\C$ of type $A_{n-1}$ will correspond to $m$-diagonals in
$\P$. Here an $m$-diagonal is a diagonal with the property that
it divides $\P$ into an $(mi +2)$-gon (for some positive integer $i$), and its complement, which is then 
an $(m(n-i)+2)$-gon.  

The actual reconstruction of the cluster category from this data, is
done in three steps:
\begin{itemize}
\item[-] construct a quiver $\Gamma$ which is isomorphic,
  as a  {\em
stable translation quiver}, to the AR-quiver of the cluster category, then
\item[-] take the {\em mesh category} of $\Gamma$, and
\item[-] take the additive category generated by the mesh category.
\end{itemize}

We shall first explain these notions, and then see how $\Gamma$ is constructed.
The AR-quiver $\Delta$ of a cluster category is an example of a stable
translation quiver. The AR-translation gives a bijective map $\tau
\colon \Delta_0 \to \Delta_0$ with the following property: 
given any two vertices $x,y$, the number of arrows $x \to y$ equals
the number of arrows $\tau y \to x$.

A (locally finite) quiver $\Gamma$ without loops, such that a
translation-function $\tau_{\Gamma}$ with the same property as $\tau$ above exists, is called
a {\em stable translation quiver}.

Given a stable translation quiver $\Gamma$ with translation function $\tau_{\Gamma}$, one can define a
{\em mesh category} $M_{(\Gamma, \tau_{\Gamma})}$. The objects in this category
are the vertices of $\Gamma$, and these are then the indecomposable
objects in the additive category generated by $M_{(\Gamma, \tau)}$.
Here we only consider quivers $\Gamma$ without multiple arrows.
In this case, the maps in $M_{(\Gamma, \tau)}$ are all linear combination of paths
modulo a certain ideal $I$ generated by the {\em mesh relations}.
For every vertex $v$ there is one mesh relation, which is constructed as
follows. Let $\{b_i \colon v_i \to v\}$ be all arrows ending in $v$, and let
$a_i \colon \tau v \to v_i$ be the arrow corresponding to $b_i$.
Then the sum $\sum b_i a_i$ is the mesh relation for $v$.

We now describe how to get a stable translation quiver $\Gamma$ from the
$(mn+2)$-gon.
Label the vertices of the polygon $1, \dots, mn+2$ (in a clockwise oriented cycle), and let $(i,j)$ denote 
an $m$-diagonal between the vertices $i$ and $j$.
We now construct a finite quiver $\Gamma$, by letting 
the vertices correspond to the $m$-diagonals. We denote by $(i,j) =
(j,i)$ the vertex corresponding to the diagonal between $i$ and $j$.
We draw an arrow $(i,j) \to (i,j+m)$, if  $(i,j+m)$ is an
$m$-diagonal,
and an arrow $(i,j) \to (i+m,j)$, if  $(i+m,j)$ is an
$m$-diagonal.
In addition, we define a translation $\tau_{\Gamma}$ by mapping $(i,j)$
to $(i- m,j-m)$.

\begin{theorem}
\cite{bm2} The $m$-cluster category of type $A_{n-1}$ is equivalent to the
additive category of the mesh category $M_{(\Gamma, \tau_{\Gamma})}$,
where $\Gamma$ is the constructed from the $(mn+2)$-gon as above.
\end{theorem}
 
\subsection{Example}
Let $m=2$ and $n=5$, and consider the $12$-gon.
It gives rise to the following stable translation quiver, which is
easily seen to be isomorphic to the AR-quiver of the $m$-cluster
category of type $A_4$ from example \ref{arex}.
 
\def \objectstyle{\scriptstyle}
\def \labelstyle{\scriptstyle}
$$
\xymatrix@!0@C0.83cm@R0.7cm{ 
 & (3,6) \ar[dr] & & (5,8) \ar[dr] & & (7,10) \ar[dr]
& & (9,12) \ar[dr]  & &  (11,2) \ar[dr]
& & (1,4) \ar[dr] & & (3,6) \\
 (1,6) \ar[dr] \ar[ur]  & & (3,8)
\ar[dr] \ar[ur]  & & (5,10) \ar[dr] \ar[ur]   & & (7,12) \ar[dr]
\ar[ur]   & & (9,2) \ar[dr] \ar[ur]  & &
(11,4)  \ar[dr] \ar[ur]  & & (1,6) \ar[dr] \ar[ur]   &  \\
& (1,8) \ar[ur]  \ar[dr]  & & (3,10) \ar[ur]  \ar[dr]  & & (5,12)
\ar[ur]  \ar[dr]   &  & (7,2) \ar[ur]  \ar[dr]  & &  (9,4) \ar[ur]  \ar[dr]  
& & (11,6) \ar[ur]  \ar[dr]  & & (1,8) \\
 (11,8) \ar[ur]  & & (1,10)
\ar[ur]  & & (3,12)  \ar[ur]   & & (5,2) 
\ar[ur]   & & (7,4)  \ar[ur]  & &
(9,6)  \ar[ur]  & &  (11,8)  \ar[ur]  &  
}
$$

\subsection{Interpretation of tilting objects and exchange}

The construction described above also has an additional important feature.
The correspondence between indecomposable objects in the $m$-cluster
category of type $A_{n-1}$ and the category of diagonals of $\P_{mn+2}$ is defined such that 
two indecomposable objects $X,Y$ in $\C$ are $\Ext$-compatible 
if and only if the $m$-diagonals corresponding to $X$ and $Y$ do
not cross. The maximal sets of non-crossing $m$-diagonals in $\P$ are called $(m+2)$-angulations.
They always have $n-1$ elements and correspond to 
the tilting objects in $\C_{A_{n-1}}$.
If we remove an $m$-diagonal in an $(m+2)$-angulation, we can replace it with $m$ different $m$-diagonals,
to obtain $m$ different $(m+2)$-angulations. This corresponds to replacing one indecomposable summand $T_k$ in
a tilting object $T$ with one of the $m$ complements of $T/T_k$ different than $T_k$.

\subsection{Example}

\begin{figure}[htp]
\begin{center}
\includegraphics[width=5cm]{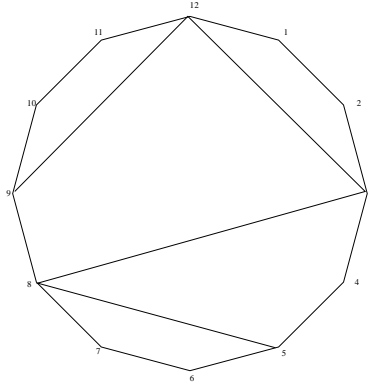}
\end{center}\caption{\label{figA} The $4$-angulation corresponding to $B \amalg X$}
\end{figure}

\sloppy The tilting object $B  \amalg  X$ of example \ref{exA}, corresponds to a $4$-angulation of a $12$-gon as in figure 
\ref{figA}. 

\begin{figure}
\begin{center}
\includegraphics[width=8cm]{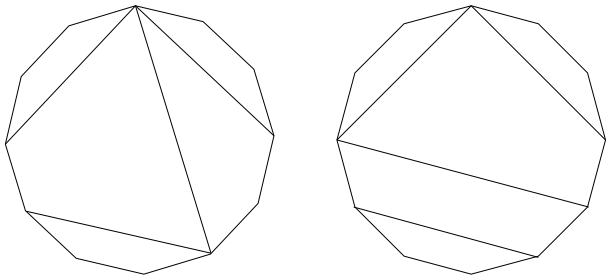}
\end{center}\caption{\label{figB} 
The $4$-angulations corresponding to $B \amalg P_2[1]$ (left) and  $B \amalg I_3[1]$ (right)}
\end{figure}

If we remove the $2$-diagonal corresponding to $X$ in this $12$-angulation, we can replace it with $m=2$ different $2$-diagonals,
and obtain the two $4$-angulations of figure \ref{figB}. These correspond to the 
tilting objects $B \amalg P_2[1]$ and $B \amalg I_3[1]$.

\section{The simplicial complex of $m$-clusters}

An {\em (abstract) simplicial complex} is a nonempty family $\Delta$ of finite
subsets of a fixed universal set,  with the property that if $X$ is in
$\Delta$, then also every subset  $Y \subset X$ is in $\Delta$.

An $m$-cluster category $\C= \C_H$ gives in a canonical way rise to a simplicial
complex $\Delta(\C)$: take the set of isomorphism classes of indecomposables in $\C$ as the universal
set,  and let $\Delta(\C)$ consist of the subsets $X$ with the property
that the elements in $X$ are $\Ext$-compatible.

Consider now the case where $\C= \C_H$ is the $m$-cluster category of $H= kQ$,
and $Q$ is a {\em Dynkin quiver}. Corresponding to the underlying
graph of $Q$ there is a finite root system $\Phi$. 

Starting with a finite root system and a positive integer $m$, Fomin and Reading \cite{fr} have defined another
simplicial complex, the {\em $m$-cluster complex}, and
one of the original motivations of studying tilting theory in
$m$-cluster categories, was to compare their simplicial complex to 
$\Delta(\C)$. This is done independently by Thomas
\cite{t} and Zhu \cite{z}. Zhu also dealt with non-simply laced Dynkin
graphs and their corresponding root systems.
The $m$-cluster complexes naturally generalises the
$1$-cluster complexes, which play a crucial role in the study of cluster
algebras \cite{fz}.

For a finite root system $\Phi$, Fomin and Reading consider the set
$\Phi_{\geq -1}^{m}$ of
{\em coloured almost positive roots}. This set consists of $m$ copies
of the positive roots, and one set of copies of the negative simple
roots. This is the universal set for the $m$-cluster complex.
Then they define a notion of compatibility of elements in this
set. This is combinatorially defined, and we leave out the details
here,
but refer instead to \cite[Section 2]{fr}.
The $m$-cluster complex consists of all sets of 
compatible elements in $\Phi_{\geq -1}^{m}$. 
Fomin and Reading show that $m$-cluster complexes satisfy some
nice conditions.

\begin{theorem}\label{fominreading}
\cite{fr} Consider a root system with $n$ simple positive roots, or equivalently
a Dynkin graph with $n$ vertices.
\begin{itemize}
\item[(a)]
All facets (inclusion-maximal sets) in the
$m$-cluster complex have cardinality $n$.
\item[(b)]
Each set in the $m$-cluster complex of cardinality $n-1$ is 
a subset of exactly $n+1$ facets.
\end{itemize}
\end{theorem}

For a given Dynkin quiver $Q$, it is well known that the set of
indecomposable $H= kQ$-modules is in bijection with the set of positive
roots of the corresponding root system. Hence, it is clear that the
indecomposable objects $\operatorname{ind} \C_H$ in the cluster
category $\C_H$ are in bijection with the set $\Phi_{\geq -1}^{m}$ of
coloured almost positive roots.

Now assume $Q$ has alternating orientation, i.e. each vertex is either
a sink or a source.
In this case Thomas \cite{t} and Zhu \cite{z} define a bijection $W$ between 
these two sets in such a way that $\Ext$-compatible objects in the
cluster category are mapped to compatible elements in
$\Phi_{\geq-1}^{m}$ . 
Hence they obtain the following.

\begin{theorem}\label{thomaszhu}
Using the bijection $W$ to identify the set of indecomposable objects in the cluster
category $\C_H$ with the set $\Phi_{\geq -1}^{m}$, the $m$-cluster
complex coincides with $\Delta(\C)$.
\end{theorem}

Let $M_{\alpha}$ be the
indecomposable $H=kQ$-module corresponding to the positive root $\alpha$.
The bijection map $W$ basically extends in a canonical way this correspondence to a correspondence
between the indecomposables in $\C$ of the form
$M[i]$, for $0 \leq i \leq m-1 $ and the $m$ copies of the positive
roots. The indecomposables $P[m] = I[-1]$ are identified with
the negative simple roots. 
 
Using this, \cite{t,z} give a conceptual, and type-free proof of
theorem \ref{fominreading}, by combining \ref{thomaszhu} with the
results in section \ref{tilt}.
Here we should note that results needed  concerning the number of
direct summands for tilting objects and the number of complements, were
proved in \cite{t,z} in the Dynkin case.

\section{Mutation of coloured quivers}

We will now discuss another combinatorial approach to $m$-cluster
categories, motivated from the fact that tilting and exchange in $1$-cluster categories
gives a categorical model for Fomin-Zelevinsky quiver mutation. We
will first recall the notion of quiver mutation.

\subsection{Fomin-Zelevinsky quiver mutation}

Let $Q = (q_{ij})$ be a finite quiver with vertices $1, \dots, n$,
with $q_{ij}$ arrows from $i$ to $j$, and with no
loops or oriented 2-cycles (parallel underlying edges with opposite directions). For a fixed  vertex $v$, we get a new
quiver $\mu_v(Q)$, also without loops or oriented two-cycles. This
operation, called {\em quiver mutation in $v$}, can be described in various
ways. Having the generalisation to $m>1$ in mind, we choose the following formulation.

\begin{itemize}
\item[-] For each pair of arrows $i \to v \to j$ in $Q$, add an arrow $i \to j$.
\item[-] If, between some pairs of vertices, there appear parallel
  underlying edges with 
  opposite directions (oriented 2-cycles), remove the same number of arrows in each
  direction, until there are no oriented 2-cycles.
\item[-] Reverse all arrows starting in or ending in $v$.
\end{itemize}

It is straightforward to check that this operation satisfies
$\mu_v(\mu_v(Q)) = Q$. It is also straightforward to verify that the quiver $\mu_v(Q)
=(\widetilde{q_{ij}} )$ is 
determined by the following formula, which is a reformulation of the
FZ-mutation formula.

\begin{equation}\label{fzfor}
\widetilde{q_{ij}} = \begin{cases} q_{ji} & \text{ if $v=i$ or $v=j$} \\
\operatorname{max} \{ 0, q_{ij} -q_{ji} + q_{iv}q_{vj} - q_{jv}q_{vi} \}
& \text{ if $i \neq v \neq j$}  
\end{cases}
\end{equation}

For a tilting object $T$ in a cluster category $\C$, we can consider 
the endomorphism-algebra $\End_{\C}(T)$. This is again a finite
dimensional basic $k$-algebra, and therefore isomorphic to a factor algebra
of a path algebra of a finite quiver $Q_T$ (the Gabriel quiver of $T$).
 
Consider now a $1$-cluster category,
let $T= B \amalg M$ and $T'=B \amalg M^{\ast}$ be two tilting objects,
and let $Q_T$ and $Q_{T^{\ast}}$ be their respective Gabriel-quivers.
The main result of \cite{bmr} is that  
\begin{equation}\label{mut} Q_{T^{\ast}} = \mu_v(Q_T),
\end{equation}
where $v$ corresponds to the indecomposable object $M$.
This can be considered a categorification of FZ-quiver mutation.

It is natural to ask for a generalisation of the above to the case $m>1$.
We give an example to show that there can be no direct
generalisation in terms of the Gabriel quiver of $T$.

\subsection{Example}\label{a2ex}
Consider the $3$-cluster category of type $A_2$. Let $P_1$ be the
simple projective, and $P_2$ be the indecomposable projective of
length $2$, and $I_2$ the simple injective. 
Then the AR-quiver of the $3$-cluster category has 11 vertices.
\bigskip
$$
\xymatrix@!0@C0.83cm@R0.7cm{ 
& P_2 \ar[dr] & & P_1[1]  \ar[dr] & & I_2[1]  \ar[dr] & & P_2[2]
\ar[dr] & & P_1[3]  \ar[dr]  & & P_1  \ar[dr] & \\
P_1 \ar[ur] & & I_2 \ar[ur] & & P_2[1] \ar[ur] & & P_1[2] \ar[ur] & &
I_2[2] \ar[ur] & & P_2[3] \ar[ur] & & P_2 
}  
$$
\bigskip
Consider the almost complete tilting object $P_2[2]$, and the four completions 
\begin{eqnarray*} & T_a  = P_2[2]   \amalg P_1 \text{,  \  \     }  & T_b
   = P_2[2]  \amalg
P_1[1] \\
  & T_c  = P_2[2]   \amalg P_1[2] \text{\ \ \ and \ \     } & T_d  = P_2[2]
 \amalg I_2[2]
\end{eqnarray*}

The following picture describes the Gabriel quivers of the
endomorphism rings of these tilting objects, with the direction of exchange indicated by the broken arrows.
$$
\xymatrix@!0@C1.00cm@R0.5cm{
T_a \colon  & \cdot & & \cdot & &  \ar@{.>}[r] & & & \cdot
& & \cdot & \colon T_b \\
& & & & & & & & & \ar@{.>}[dd]  & & \\
& &    & & & &  & & &  & & \\
& & \ar@{.>}[uu]& && & & & & & & \\
T_d \colon & \cdot  & & \cdot \ar[ll] & & & \ar@{.>}[l]   &
& \cdot \ar[rr]& & \cdot  &   \colon T_c 
}
$$
From this it is clear that more information than the Gabriel quiver of a tilting
object $T$ is needed, in order to generalise formula (\ref{mut}). 

\subsection{Coloured quivers and mutation}

It turns out that instead of Gabriel quivers, we can now deal with
coloured quivers. 

An $m$-coloured multi-quiver $Q$, consists of vertices $1, \dots, n$
and coloured arrows $i \overset{(c)}{\to} j$, where $c$ is in $\{0,1,
\dots, m\}$. We let $q_{ij}^{(c)}$ denote the number of arrows from
$i$ to $j$ of colour $(c)$.

Coloured quiver mutation was introduced in \cite{bt}.
Given a vertex $v$ in an $m$-coloured quiver $Q$, define
a new coloured quiver $\mu_v(Q)$ by modifying $Q$ as follows.

\begin{itemize}
\item[-] For each pair of arrows 
$$\xymatrix{
i \ar[r]^{(c)} & v \ar[r]^{(0)} & j 
} 
$$
with $c$ in $0, 1, \dots m$, add two arrows: one arrow of colour $(c)$ from $i$ to
$j$ and one arrow of colour $(m-c)$ from $j$ to $i$.
\item[-]  If, for some pairs of vertices, there appear parallel arrows 
with different colours from $i$ to $j$, remove the same number of arrows of each
colour. 
\item[-] Change the colour of all arrows ending in $v$, by adding
  one. Change the colour of all arrows starting in $v$, by subtracting one.
\end{itemize}

Alternatively one can describe coloured mutation via a formula which
is a generalised version of formula (\ref{fzfor}).
If $Q= (q_{ij})$ is an $m$-coloured quiver, then
$Q' = \mu_v{Q} = (\widetilde{q_{ij}})$  is given by \footnote{Note that in
  \cite{bt}, there is an unfortunate typo in the formula: the two
  first cases are mixed up.}:
$$
{\scriptstyle\widetilde{q}_{ij}^{(c)} = 
\begin{cases}  q_{ij}^{(c+1)} & \text{  if $v =i$} \\
		 q_{ij}^{(c-1)} &\text{  if $v=j$} \\
		 \max \{0, q_{ij}^{(c)} - \sum_{t \neq c} q_{ij}^{(t)} + (q_{iv}^{(c)} - q_{iv}^{(c-1)}) q_{vj}^{(0)} 
		 + q_{iv}^{(m)} (q_{vj}^{(c)}  -q_{vj}^{(c+1)}) \} & \text{  else} 
                 \end{cases}
}
$$

\subsection{The coloured quiver of a tilting object}
Let $m \geq 1$ be an integer, and $\C$  an $m$-cluster category.
We want to assign to each tilting object $T = \amalg_{i=1}^n T_i$ in $\C$ a coloured quiver
$Q_T = (q_{ij}^{(c)})$ with $n$ vertices corresponding to the 
indecomposable direct summands in $T$. To determine the coloured
arrows, we use the exchange triangles 
(\ref{exch}): we let $q_{ij}^{(c)}$ be the multiplicity of $T_j$ as a
direct summand in $B_i^{(c)}$. Note that the $0$-coloured arrows
are indeed the arrows of the Gabriel quiver of $T$. 

Not all coloured quivers can be obtained as $Q_T$ for a 
tilting object $T$. By definition, there are no loops (of any colour) in $Q_T$,
that is: $q_{ii}^{(c)} = 0$ for $i$ and all $c$. Also, one can prove
that $Q_T$ is locally monochromatic: for fixed vertices $i,j$ there are only
arrows of one colour from $i$ to $j$. One can also prove that 
$q_{ij}^{(c)} = q_{ji}^{(m-c)} $, that is: for each arrow of colour $c$, there is an
arrow in the opposite direction with  colour $m-c$. 
There are also more known restrictions, see \cite[Prop. 5.1]{bt}.

It is an interesting open problem to find a set of properties that
characterises the coloured quivers of type $Q_T$ among all coloured
quivers.

One can now generalise the result in \cite{bmr} to coloured quivers of
tilting objects in higher cluster categories.

\begin{theorem}\label{colmut}
Let $T= \amalg_{i=1}^n T_i$ and $T' = T/T_j \amalg T_j^{(1)}$ be
tilting objects in an $m$-cluster category $\C$, such that there is
an exchange triangle 
\begin{equation}  T_j \rightarrow B_j^{(0)}
\rightarrow T_j^{(1)} \rightarrow.  \end{equation} 
Then $Q_{T'} = \mu_j(Q_T)$.
\end{theorem}

\subsection{Example}
Revisiting example \ref{a2ex}, we now consider instead the coloured
quivers, and their mutations. Note that we always mutate in the
leftmost vertex.
$$
\xymatrix@!0@C1.00cm@R0.5cm{
T_a\colon & \cdot \ar@<0.3ex>_{(1)}[rr] & & \cdot
\ar@<0.3ex>_{(2)}[ll]  & & \ar@{.>}[r] & & &  \cdot \ar@<0.3ex>_{(2)}[rr] & & \cdot
\ar@<0.3ex>_{(1)}[ll] &  \colon T_b  \\
& & & & & &  & & & \ar@{.>}[dd]   & &  \\
& &    & & & &   & & & & & \\
&  & \ar@{.>}[uu]&  & & & & & & & & \\
T_d \colon & \cdot \ar@<0.3ex>_{(0)}[rr] & & \cdot
\ar@<0.3ex>_{(3)}[ll]  & & &  \ar@{.>}[l]   & & \cdot \ar@<0.3ex>_{(3)}[rr] & & \cdot
\ar@<0.3ex>_{(0)}[ll]   & \colon T_c  
}
$$

\subsection{Example}
We consider again the case $m=2$, with the quiver $Q$ of type $A_4$ as
in example \ref{ex1}.

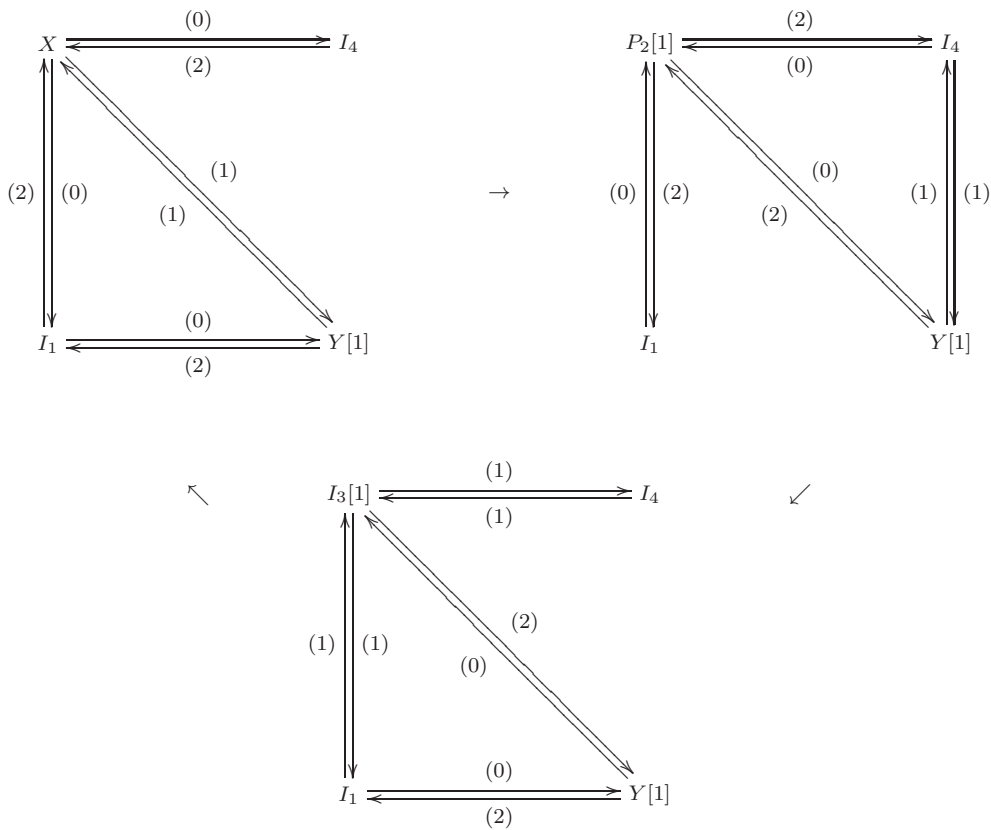
\begin{figure}

\def \labelstyle{\scriptstyle}
$$
\xymatrix@!0@C2.0cm@R2.0cm{ 
X \ar@<0.3ex>^{(0)}[rr] \ar@<0.3ex>^{(0)}[dd]  \ar@<0.3ex>^{(1)}[ddrr]
& & I_4
\ar@<0.3ex>^{(2)}[ll]  & &  P_2[1] \ar@<0.3ex>^{(2)}[rr]
\ar@<0.3ex>^{(2)}[dd]  \ar@<0.3ex>^{(0)}[ddrr] & & I_4
\ar@<0.3ex>^{(0)}[ll] \ar@<0.3ex>^{(1)}[dd] & \\
& & &  \to &  & &  \\
I_1 \ar@<0.3ex>^{(0)}[rr] \ar@<0.3ex>^{(2)}[uu]  & & Y[1] \ar@<0.3ex>^{(2)}[ll]
\ar@<0.3ex>^{(1)}[uull] &  & I_1 \ar@<0.3ex>^{(0)}[uu]  & & Y[1] \ar@<0.3ex>^{(1)}[uu] \ar@<0.3ex>^{(2)}[uull] \\ 
& \nwarrow & I_3[1] \ar@<0.3ex>^{(1)}[rr] \ar@<0.3ex>^{(1)}[dd]
\ar@<0.3ex>^{(2)}[ddrr] &  & I_4 \ar@<0.3ex>^{(1)}[ll]  &  \swarrow &\\
& & & & & &\\
& & I_1 \ar@<0.3ex>^{(1)}[uu]  \ar@<0.3ex>^{(0)}[rr]  &  & Y[1]
\ar@<0.3ex>^{(0)}[uull] \ar@<0.3ex>^{(2)}[ll]  & &
}
$$
\caption{{Coloured mutation at the upper left vertex}\label{fig1}}
\end{figure}

The coloured quivers of the three tilting objects
$$ T  =   I_1 \amalg I_4 \amalg Y[1] \amalg X \text{, \
  \  \ }
 T'  =  I_1 \amalg I_4 \amalg Y[1] \amalg P_2[1] \text{   and   } $$  
\ 
$$ T''  =  I_1 \amalg I_4 \amalg Y[1] \amalg I_3[1]$$
are given in Figure \ref{fig1}. Note that $Q_{T'}$ is given by coloured mutation
of  $Q_{T}$ at the vertex corresponding to $X$, that 
 $Q_{T''}$ is given by coloured mutation
of  $Q_{T'}$ at the vertex corresponding to $P_2[1]$, and that
$Q_{T}$ is given by coloured mutation
of  $Q_{T''}$ at the vertex corresponding to $I_3[1]$.

\subsection{Finiteness of the mutation class}

Let $Q$ be an acyclic quiver.
We can view this as an $m$-coloured quiver, by regarding each arrow $\alpha$ in $Q$ as
an arrow of colour $(0)$, and then adding an arrow of colour
$(m)$ in opposite direction to $\alpha$.

Torkildsen \cite{tor1} has proved the following, generalising
a similar statement of \cite{br} for $m=1$. 

\begin{theorem}
\cite{tor1} The coloured mutation class of a connected acyclic quiver $Q$ is finite if and
only if $Q$ is either of Dynkin or extended Dynkin type, or has at
most two vertices.
\end{theorem}

In Dynkin type $A$,  Torkildsen \cite{tor2} has also found a formula for the number of elements in
the mutation class, using a connection to the classical 
cell-growth
problem \cite{hpr}. Fomin and Reading \cite{fr} have shown that number of
$m$-clusters (in the Dynkin case) is given by the Fuss-Catalan numbers.

\subsection{$m$-cluster tilted algebras}

Coloured quiver mutation gives some information on the $m$-cluster-tilted
algebras, i.e.  
algebras of the form $\End_{\C}(T)$ for $T$ a cluster-tilting object
in an $m$-cluster category.

Using that any tilting object can be reached from
any other tilting object by a sequence of exchanges \cite{zz}, one
obtains the following as a consequence of Theorem \ref{colmut}.

\begin{theorem}
 \cite{bt} Let $\C= \C_{kQ}$ for an
acyclic quiver $Q$. Then the Gabriel quivers of all $m$-cluster tilted
algebras are obtained by iterated coloured mutation of $Q$.
\end{theorem}

\section{Other aspects and generalisations}

In this survey, the main focus is on the combinatorial aspects of
higher cluster categories. In this concluding section, we give some
links to other aspects and generalisations, leaving out all details.

\subsection{Calabi-Yau triangulated categories}\label{cy}

Consider a triangulated category $\C$ with split idempotents and with 
suspension functor $\Sigma$. Assume in addition that all $\Hom$-spaces
of $\C$ are finite dimensional over the algebraically closed field $k$, and that $\C$ admits a Serre functor
$\nu$, i.e. there is a bifunctorial isomorphism $$\Hom_{\C}(X, \nu Y) \simeq
D\Hom_{\C}(Y,X).$$
If, in addition, there is an isomorphism $\Sigma^{m+1} \simeq \nu$, then $\C$ is said
to be Calabi-Yau of CY-dimension $m+1$ (for short $m+1$-Calabi-Yau).
Note that the $m$-cluster
category satisfies all these properties with $\nu = \tau [1]$.

Rigid objects and tilting objects may now be defined exactly as in the 
case of $m$-cluster categories. In fact, one does not need to restrict
to objects. In \cite{kr1}, a {\em (cluster) tilting subcategory}
in a $m+1$-Calabi-Yau category is defined as a 
$k$-linear functorially finite subcategory $\T$ of $\C$, satisfying

\begin{itemize}
\item[-] $\Ext^i(T,T') = 0$ for all $T,T'$ in $\T$ and all $0<i<m$, and
\item[-] if $X \in \C$ satisfies $\Ext^i(T,X) = 0$ for all $T$ in $\T$
  and all $0<i<m$, then $X$ belongs to $\T$.
\end{itemize}

Note that the additive closure $\add T$ of a tilting object $T$ in an
$m$-cluster category clearly satisfies this.
Keller and Reiten \cite{kr2}, showed that one can characterise
$m$-cluster categories as exactly those
$m+1$-Calabi-Yau categories with an object $T$, such that 
\begin{itemize}
\item[-] $\add T$ is a cluster tilting subcategory
\item[-] $\Hom(T,\Sigma^i T) = 0$ for $i= -m, \dots , -1$, and
\item[-] $\End(T)$ is a hereditary algebra.
\end{itemize}

\subsection{Generalised higher cluster categories}

Amiot gave in \cite{amiot} a more general definition of cluster
categories in the case $m=1$. Starting with a finite dimensional
algebra $A$ of global dimension at most $2$, she constructs a certain
triangulated category $\C_A$, which is equivalent to the ordinary
cluster category in case $A$ is hereditary. This category $\C_A$ is
in general not $\Hom$-finite. But, if $A$ satisfies certain additional
conditions, then $\C_A$ is $\Hom$-finite, and in this case
$\C_A$ is $2$-Calabi-Yau and 
$A$ is a tilting object in $\C_A$.

In a very recent paper Lingyan Guo \cite{guo} generalises this construction to 
 $m >1$. More precisely; for finite dimensional algebra $A$ of finite
global dimension $m$, assume that the functor $\Tor^A_m (- , DA)$ is nilpotent.  
In this setting she
constructs a $\Hom$-finite triangulated category $\C_A^{(m-1)}$, which is
$m$-Calabi-Yau, and such that $A$ is an
$m-1$-cluster tilting
object in  $\C_A^{(m-1)}$.

In addition, both in \cite{amiot} and \cite{guo}, generalised (higher) cluster
categories are also considered in the setting of quivers with
(super-)potentials, see \cite{dwz}.

%



\end{document}